\documentclass[10pt,a4paper]{article}

\usepackage[latin2]{inputenc}
\usepackage{amsmath}
\usepackage{amsfonts}
\usepackage{amssymb}
\usepackage{makeidx}

\newtheorem{theorem}{Theorem}
\newtheorem{conjecture}{Conjecture}

\def\stop{\mbox{\footnotesize {\vrule width 6pt height 6pt}}}
\def\nfrac#1#2{\mbox{\footnotesize $\displaystyle\frac{#1}{#2}$}}

\begin{document}

$\,$

\bigskip

\centerline{\large \bf A PROOF OF TWO CONJECTURES OF CHAO$\!\;$-PING CHEN}

\smallskip

\centerline{\large \bf FOR INVERSE TRIGONOMETRIC FUNCTIONS}

\medskip
\begin{center}
{\sc Branko Male\v sevi\' c}${\;}^{1)}$, {\sc Bojan Banjac}${\;}^{1,2)}$, {\sc Ivana Jovovi\' c}${\;}^{1)}$
\end{center}

\begin{center}
{\small ${\;}^{1)}\,$University of Belgrade, Faculty of Electrical Engineering,   \\[0.0 ex]
Department of Applied Mathematics, Serbia}

\medskip

{\small ${\;}^{2)}\,$University of Novi Sad, Faculty of Technical Sciences,       \\[0.0 ex]
Computer Graphics Chair, Serbia}
\end{center}

\medskip
\noindent
{\small \textbf{Abstract.}
In this paper we prove two conjectures stated by Chao\,-Ping Chen in [Int.~Trans.~Spec.~Funct. 23:12 (2012), 865--873],
using a method for proving inequalities of mixed trigonometric polynomial functions.

\medskip
{\footnotesize Keywords: inequalities, inverse trigonometric functions, Taylor series}

\footnote{$\!\!\!\!\!\!\!\!\!\!\!\!\!\!$
\scriptsize Emails: Branko Male\v sevi\' c {\tt $<$malesevic@etf.bg.ac.rs$>$}, Bojan Banjac {\tt $<$bojan.banjac@uns.ac.rs$>$},
Ivana Jovovi\' c {\tt $<$ivana@etf.bg.ac.rs$>$}}

\vspace*{-2.5 mm}

{\small \tt MSC:26D05}

\bigskip

\bigskip
\noindent
Wilker in \cite{Wilker_1989} formulated two problems.
First one was to prove that
\begin{equation}
\label{wilker_osnovno}
\left(\nfrac{\sin x}{x}\right)^{\!2}+\nfrac{\tan x}{x}>2
\end{equation}
holds for $0 < x < \nfrac{\pi}{2}$;
and second one was to find the largest constant $c$ such that
$$
\left(\nfrac{\sin x}{x}\right)^{\!2}+\nfrac{\tan x}{x}>2+c\,x^3\tan x,
$$
for $0<x<\nfrac{\pi}{2}$.

\medskip
\noindent
Sumner, Jagers, Vowe, and Anglesio in \cite{Sumner_Jagers_Vowe_Anglesio_1991} gave an improvement
of the inequality (\ref{wilker_osnovno}) in the form
$$
2+\left(\nfrac{2}{\pi}\right)^{\!4}x^3\tan x
<
\left(\nfrac{\sin x}{x} \right)^{\!2}+\nfrac{\tan x}{x}<2+\nfrac{8}{45}x^3\tan x,
$$
where the constants $\left(\nfrac{2}{\pi}\right)^{\!4}$ and $\nfrac{8}{45}$
were the best possible.
Huygens in \cite{Huygens_1888} presented the inequality
\begin{equation}
\label{Huygens_osnovno}
2\;\nfrac{\sin x}{x} + \nfrac{\tan x}{x} > 3,
\end{equation}
for $0<|x|<\nfrac{\pi}{2}$.

\medskip
\noindent
Neuman, and Sandor in \cite{Neuman_Sandor_2010} established the relation between inequalities (\ref{wilker_osnovno}) and (\ref{Huygens_osnovno}).
The relevant papers on the topic are also \cite{Zhu_2005}, \cite{Zhang_Zhu_2007}, \cite{Sun_Zhu_2011}, \cite{Sun_Zhu_2011b} and \cite{Debnath_Mortici_Zhu_2015}.
The inverse trigonometric and inverse hyperbolic versions of Wilker and Huygens's inequalities were considered
in \cite{Zhu_2007}, \cite{Neuman_Sandor_2010}, \cite{C.-P. Chen and W.-S. Cheung 2012} and \cite{C.-P. Chen and W.-S. Cheung 2012 (JIA)}.
Recently, the analogue inequalities for the generalized trigonometric functions \cite{Neuman_2014 (AMC)} and different special functions \cite{Neuman_2014}~and~\cite{Neuman_2015}
have been taken into consideration.

\medskip
\noindent
Chao\,-Ping Chen in \cite{Chao-Ping Chen} proved the following two theorems and proposed two open problems.
\begin{theorem}
If $0< x < 1$, then
$$
2 + \nfrac{17}{45} \, x^3 \arctan x
<
\left(\nfrac{\arcsin x}{x}\right)^{\!2}
+
\nfrac{\arctan x}{x},
$$
where the constant $\displaystyle\nfrac{17}{45}$ is the best possible.
\end{theorem}
Considering the previous theorem it was natural to ask what is
the best possible constat $c$ such that
$$
\left(\nfrac{\arcsin x}{x}\right)^{\!2}
+
\nfrac{\arctan x}{x}
<
2 + c \,x^3 \arctan x
$$
holds, for $0\!<\!x\!<\!1$. The choice of constant $\left(\pi^2\!+\!\pi\!-\!8\right)/\pi$ is somehow motivated,
since it is the limit at $\nfrac{\pi}{2}$ of the function $x \mapsto \left(\left(\arcsin x / x\right)^{2}
\!+\left(\arctan x/x\!-\!2\right)\right) / \left(x^3 \arctan x\right)$. Therefore, Chao\,-Ping Chen
in \cite{Chao-Ping Chen} stated the following conjecture.
\begin{conjecture}
\label{first conjecture}
If $0 < x < 1$, then
\begin{equation}
\label{first inequality}
\left(\nfrac{\arcsin x}{x}\right)^{\!2}
+
\nfrac{\arctan x}{x}
<
2 + \nfrac{\pi^2+\pi-8}{\pi} \, x^3 \arctan x,
\end{equation}
where the constant $\displaystyle\nfrac{\pi^2+\pi-8}{\pi}$ is the best possible.
\end{conjecture}

\smallskip
\noindent
In the paper \cite{Chao-Ping Chen} one can also find the following theorem.
\begin{theorem}
\label{second conjecture}
If $0< x < 1$, then
$$
3 + \nfrac{7}{20} \, x^3 \arctan x
<
2\left(\nfrac{\arcsin x}{x}\right)
+
\nfrac{\arctan x}{x},
$$
where the constant $\displaystyle\nfrac{7}{20}$ is the best possible.
\end{theorem}
And so, there is a matched conjecture.
\begin{conjecture}
If $0< x < 1$, then
\begin{equation}
\label{second inequality}
2\left(\nfrac{\arcsin x}{x}\right)
+
\nfrac{\arctan x}{x}
<
3 + \nfrac{5\,\pi-12}{\pi} \, x^3 \arctan x,
\end{equation}
where the constant $\displaystyle\nfrac{5\,\pi-12}{\pi}$ is the best possible.
\end{conjecture}
The proofs of the previous two theorems are based on the usage of the appropriate infinite power series.
In the proofs of the stated conjectures a method from \cite{Malesevic_Makragic_2015}
will be used and it is based on the usage of the appropriate approximations of some mixed trigonometric polynomials with finite Taylor series.
This method presents continuation of Mortici's method from \cite{Mortici_2011}.
The method is also applied on inequalities closely related to presented ones,
see \cite{Banjac_Makragic_Malesevic_2015} and \cite{Nenezic_Malesevic_Mortici_2015}.

\bigskip
\noindent
We follow the notation used in \cite{Malesevic_Makragic_2015}.
Let $\varphi\!:\![a,b]\!\longrightarrow\!\mathbb{R}$ be a differentiable function on a segment $[a,b]$
and differentiable on a right at $x\!=\!a$ an arbitrary number of times.
Denote by $T^{\varphi,\, a}_{m}(x)$ the Taylor polynomial of the order $m$ of the function $\varphi$ in the point $x\!=\!a$.
If there is some $\eta\!>\!0$ such that \mbox{$T^{\varphi, \,a}_{m}(x) \geq \varphi(x)$} holds for $x \!\in\! (a,a+\eta) \!\subset\! [a,b]$,
then we define $\overline{T}^{\,\varphi,\, a}_{m}(x)\!=\!T^{\,\varphi,\, a}_{m}(x)$,
and $\overline{T}^{\,\varphi,\, a}_{m}(x)$ presents an upward approximation of the order $m$ of the function $\varphi$
in the right neighborhood $(a,a+\eta)$ of the point $a$.
Analogously, if there is some $\eta>0$ such that \mbox{$T^{\varphi, a}_{m}(x) \leq \varphi(x)$} holds for $x \!\in\! (a,a+\eta) \!\subset\! [a,b]$,
then we define $\underline{T}^{\varphi,\, a}_{\,m}(x)\!=\!T^{\,\varphi,\, a}_{m}(x)$,
and $\underline{T}^{\varphi,\, a}_{\,m}(x)$ presents a downward approximation of the order $m$ of the function $\varphi$
in the right neighborhood $(a,a+\eta)$ of the point $a$.
In the same manner, it is possible to define upward and downward approximations in the left neighborhood of a point.

\section{Proof of the Conjecture 1} 

Let us first observe the inequality (\ref{first inequality}) of the Conjecture \ref{first conjecture} written in the form
\begin{equation}
\label{startinequality0}
2+\nfrac{\pi^2+\pi-8}{\pi}\,x^3\arctan x - \nfrac{\arctan x}{x}-\left(\nfrac{\arcsin x}{x}\right)^{\!2}
>
0,
\end{equation}
for $x \in (0,1)$.
Substituting $x=\sin t$ into (\ref{startinequality0}), for $t \in \left(0,\nfrac{\pi}{2}\right)$,
we obtain
\begin{equation}
\label{startinequality1}
2+
\nfrac{(\pi^2+\pi-8)\sin^4 t - \pi}{\pi \sin t} \; \arctan(\sin t) - \nfrac{t^2}{\sin^2 t} > 0.
\end{equation}
It is enough to prove that
\begin{equation}
\label{startinequality2}
\begin{array}{rcl}
g(t)
\!\!&\!=\!&\!\!
2\,\pi \sin^2 t
+
(\pi^2+\pi-8) \sin^5\!t\;\arctan(\sin t)                                        \\[1.0 ex]
\!\!&\! \!&\!\!
-\,
\pi \, \sin t \; \arctan(\sin t)
-
\pi\,t^2
>
0,
\end{array}
\end{equation}
for $t \in \left( 0, \nfrac{\pi}{2} \right)$.
Let us notice that $t=0$ is zero of the sixth order and $t=\nfrac{\pi}{2}$ is the simple zero of the function $g$.
Furthermore, we differentiate two cases if $t \!\in\! (0, 1.1]$ or $t \!\in\! (1.1,\pi/2)$.

\medskip  
\noindent
{\bf (I)} {\boldmath $t \!\in\! (0, 1.1]$}
Let us start from the series
$
\arctan x
\!=\!
\mbox{\footnotesize $\displaystyle\sum$}_{k=0}^{\infty}{\mbox{\footnotesize $(-1)^k$}\nfrac{x^{2k+1}}{2k+1}}
$,
which holds for $x \!\in\! [-1,1]$.
Notice that
$
\underline{T}^{\mbox{\footnotesize \rm arctan}, 0}_{3+4k_1}(x)
<
\arctan x
<
\overline{T}^{\mbox{\footnotesize \rm arctan}, 0}_{1+4k_2}(x)
$
are true for $x \in (0,1]$ and $k_{1,2} \!\in\! \mathbb{N}_{0}$.
By introducing the substitution $x = \sin t$, we can conclude that
$$
\label{dobleinequality_general}
\underline{T}^{\mbox{\footnotesize \rm arctan}, 0}_{3+4k_1}(\sin t)
\!<\!
\arctan(\sin t)
\!<\!
\overline{T}^{\mbox{\footnotesize \rm arctan}, 0}_{1+4k_2}(\sin t),
$$
for $t \in \left(0,\nfrac{\pi}{2}\right)$ and $k_{1,2} \!\in\! \mathbb{N}_{0}$.
For the proof of the Conjecture \ref{first conjecture} we will use previous inequalities for $k_1=0$ and $k_2=1$,
i.e. we will only need
\begin{equation}
\label{dobleinequality_4_6}
\underline{T}^{\mbox{\footnotesize \rm arctan}, 0}_{3}(\sin t)
<
\arctan(\sin t)
<
\overline{T}^{\mbox{\footnotesize \rm arctan}, 0}_{5}(\sin t),
\end{equation}
for $t \in \left(0,\nfrac{\pi}{2}\right)$.
Since $\underline{T}^{\mbox{\footnotesize \rm arctan}, 0}_{3}(\sin t) > 0$
for $t \in \left(0,\nfrac{\pi}{2}\right)$ and $\pi^2+\pi-8 > 0$ we have
$$
\begin{array}{rcl}
g(t)
\!\!&\!>h(t)\;=\;\!&\!\!
2\,\pi \sin^2\!t
+
(\pi^2+\pi-8) \sin^5\!t \; \underline{T}^{\mbox{\footnotesize \rm arctan}, 0}_{3}(\sin t)     \\[1.0 ex]
\!\!&\! \!&\!\!
-
\pi \sin t \; \overline{T}^{\mbox{\footnotesize \rm arctan}, 0}_{5}(\sin t)
-
\pi \, t^2,
\end{array}
$$
for $t \in \left(0,1.1\right]$.
It remains to prove that
$$
\begin{array}{rcl}
h(t)
\!\!&\!=\!&\!\!
2\,\pi \sin^2\!t
+
(\pi^2+\pi-8) \sin^5\!t
\left(\sin t - \nfrac{1}{3} \sin^3\!t\right )                                   \\[1.0 ex]
\!\!&\! \!&\!\!
-\,
\pi \sin t
\left( \sin t - \nfrac{1}{3} \sin^3\!t + \nfrac{1}{5} \sin^5\!t\right )
-
\pi \, t^2
>
0,
\end{array}
$$
for $t \in \left(0,1.1\right]$. The function $h$ is a mixed trigonometric polynomial function.
For the proof of the inequality $h(t) > 0$, for $t \in \left( 0, 1.1 \right]$,
we use method from the paper \cite{Malesevic_Makragic_2015}.
Using trigonometric multiple angle formulas, we obtain
$$
\begin{array}{rcl}
h(t)
\!&\!\!=\!\!&\!
\left(
-\nfrac{\pi^2}{384}
-\nfrac{\pi}{384}
+\nfrac{1}{48}
\right)\cos 8\,t
+
\left(
-\nfrac{\pi^2}{96}
-\nfrac{\pi}{240}
+\nfrac{1}{12}
\right)\cos 6\,t                                                                \\[1.0 ex]
\!&\! \!&\!
+
\left(
\nfrac{11\,\pi^2}{96}
+\nfrac{19\,\pi}{160}
-\nfrac{11}{12}
\right)\cos 4\,t
+
\left(
-\nfrac{31\,\pi^2}{96}
-\nfrac{43\,\pi}{48}
+\nfrac{31}{12}
\right)\cos 2\,t                                                                \\[1.0 ex]
\!&\! \!&\!
-
\pi \,t^2
+
\nfrac{85\,\pi^2}{384}
+
\nfrac{301\,\pi}{384}
-
\nfrac{85}{48}.
\end{array}
$$
Inequalities from the paper \cite{Malesevic_Makragic_2015}:
$$
\overline{T}^{\mbox{\footnotesize \rm cos}, 0}_{k}(y)
>
\cos y\,\,{\big (}k \!=\! 4, 12, 16{\big )}
\;\;\;\mbox{and}\;\;\;
\cos y
>
\underline{T}^{\mbox{\footnotesize \rm cos}, 0}_{10}(y),
$$
$y \in {\big (}\,0,\sqrt{(k+3)(k+4)}\;{\big )}$, yield
$$
\begin{array}{rcccl}
h(t)
\!&\!>\!&\!
P_{16}(t)
\!&\!=\!&\!
\left(
-\nfrac{\pi^2}{384}
-\nfrac{\pi}{384}
+\nfrac{1}{48}
\right)\overline{T}^{\mbox{\footnotesize \rm cos}, 0}_{16}(8\,t)                \\[1.75 ex]
\!&\! \!&\!
\!&\! \!&\!
\!+
\left(
-\nfrac{\pi^2}{96}
-\nfrac{\pi}{240}
+\nfrac{1}{12}
\right)\overline{T}^{\mbox{\footnotesize \rm cos}, 0}_{12}(6\,t)                \\[1.75 ex]
\!&\! \!&\!
\!&\! \!&\!
\!+
\left(
\nfrac{11\,\pi^2}{96}
+\nfrac{19\,\pi}{160}
-\nfrac{11}{12}
\right)\underline{T}^{\mbox{\footnotesize \rm cos}, 0}_{10}(4\,t)               \\[1.75 ex]
\!&\! \!&\!
\!&\! \!&\!
\!+
\left(
-\nfrac{31\,\pi^2}{96}
-\nfrac{43\,\pi}{48}
+\nfrac{31}{12}
\right)\overline{T}^{\mbox{\footnotesize \rm cos}, 0}_{4}(2\,t)                 \\[1.75 ex]
\!&\! \!&\!
\!&\! \!&\!
\!-
\pi \,t^2
+
\nfrac{85\,\pi^2}{384}
+
\nfrac{301\,\pi}{384}
-
\nfrac{85}{48},
\end{array}
$$
for $t \in (0,1.1]$.
Hence we prove that
$$
\begin{array}{rcl}
P_{16}(t)
\!&\!=\!&\!
\left(
-
\nfrac{67108864\,\pi^2}{1915538625}
-
\nfrac{67108864\,\pi}{1915538625}
+
\nfrac{536870912}{1915538625}
\right) t^{16}                                                                  \\[1.75 ex]
\!&\! \!&\!
+
\left(
\nfrac{16777216\,\pi^2}{127702575}
+
\nfrac{16777216\,\pi}{127702575}
-
\nfrac{134217728}{127702575}
\right)t^{14}                                                                   \\[1.75 ex]
\!&\! \!&\!
+
\left(
-
\nfrac{945149\,\pi^2}{2245320}
-
\nfrac{11017201\,\pi}{28066500}
+
\nfrac{945149}{280665}
\right) t^{12}                                                                  \\[1.75 ex]
\!&\! \!&\!
+
\left(
\nfrac{309929\,\pi^2}{340200}
+
\nfrac{27409\,\pi}{34020}
-
\nfrac{309929}{42525}
\right) t^{10}                                                                  \\[1.75 ex]
\!&\! \!&\!
+
\left(
-
\nfrac{20129\,\pi^2}{15120}
-
\nfrac{1609\,\pi}{1512}
+
\nfrac{20129}{1890}
\right) t^8                                                                     \\[1.75 ex]
\!&\! \!&\!
+
\left(
\nfrac{1049\,\pi^2}{1080}
+
\nfrac{293\,\pi}{540}
-
\nfrac{1049}{135}
\right) t^6
>
0
\end{array}
$$
holds, for $t \in (0,1.1]$. Notice that $P_{16}(t) = \nfrac{t^6}{30648618000} P_{10}(t)$ for
$$
\begin{array}{rcl}
P_{10}(t)
\!&\!=\!&\!
{\big (}\!-\!1073741824\,\pi^2\!-\!1073741824\,\pi\!+\!8589934592{\big )}\,t^{10}  \\[1.0 ex]
\!&\! \!&\!
\!+
{\big (}4026531840\,\pi^2\!+\!4026531840\,\pi\!-\!32212254720{\big )}\,t^8         \\[1.0 ex]
\!&\! \!&\!
\!+
{\big (}\!-\!12901283850\,\pi^2\!-\!12030783492\,\pi\!+\!103210270800{\big )}\,t^6 \\[1.0 ex]
\!&\! \!&\!
\!+
{\big (}27921503610\,\pi^2\!+\!24692768100\,\pi\!-\!223372028880{\big )}\,t^4      \\[1.0 ex]
\!&\! \!&\!
\!+
{\big (}\!-\!40801986225\,\pi^2\!-\!32614832250\,\pi\!+\!326415889800{\big )}\,t^2 \\[1.0 ex]
\!&\! \!&\!
\!+
{\big (}29768889150\,\pi^2\!+\!16629713100\,\pi\!-\!238151113200{\big )}.
\end{array}
$$
Let us introduce substitution $z\!=\!t^2$, for $z \in (0, 1.21]$, and prove that $P_{5}(z) = P_{10}(\sqrt{z}\,)> 0$.
According to the Ferrari's formulas, the derivative polynomial $P_{5}'$ does not have real roots.
Since $P_{5}'(0) < 0$ we can assert that $P_{5}'(z) < 0$, for every $z \in (0, 1.21]$.
Therefore, $P_{5}$ is strictly decreasing function with unique real root $z_1 = 1.233 \ldots > 1.21$.
So we have $P_{5}(z) > 0$ for $z \in (0, 1.21]$, i.e. $P_{10}(t) > 0$ for $t \in (0,1.1]$.
Finally, we conclude
$$
g(t) > h(t) > P_{16}(t) = \nfrac{\,t^6}{30648618000} P_{10}(t) > 0,
$$
for $t \in (0,1.1]$.

\medskip  
\noindent
{\bf (II)} {\boldmath $t \!\in\! \left(1.1, \pi\!\,/\,\!2\right)$}
We transform the inequality (\ref{startinequality2}), for $t \in \left(1.1,\nfrac{\pi}{2}\right)$, to the inequality
$$
\label{startinequality3}
\begin{array}{rcl}
g_{2}(t)
=
g(\nfrac{\pi}{2} - t)
\!\!&\!=\!&\!\!
{\Big (}(\pi^2+\pi-8)\cos^4\!t-\pi{\Big )}\cos t \;\arctan(\cos t)              \\[1.0 ex]
\!\!&\! \!&\!\!
-\,\pi\,\left(\nfrac{\pi}{2}-t\right)^{\!2} + 2\,\pi \cos^2 t
>
0,
\end{array}
$$
for $t \in \left( 0, \nfrac{\pi}{2} - 1.1 \right) = {\big (} 0, 0.470 \ldots {\big )}$.
Let us notice that
$(\pi^2+\pi-8)\cos^4\!t-\pi > 0$ \hfill

\smallskip
\noindent
is true for $t \in \left( 0, \nfrac{\pi}{2} - 1.1 \right)$.
Furthermore, we consider the additional inequality
\begin{equation}
\label{auxiliaryineq}
\arctan(\cos t) \geq \nfrac{\pi}{4}-\nfrac{t}{2},
\end{equation}
for $t \in \left[0, \nfrac{\pi}{2}\right]$.
Equality holds for $t=0$ or $t=\nfrac{\pi}{2}$.
Obviously
$$
\left(\arctan(\cos t)\right)' \!=\! \nfrac{-\sin t}{\cos^2 t\!+\!1}\!<\!0
\;\;\;\mbox{and}\;\;\;
\left(\arctan(\cos t)\right)'' \!= \! \nfrac{(\cos^2t\!-\!3) \cos t}{\;(\cos^2 t\!+\!1)^2}\!<\!0,
$$
for $t \in \left[0, \nfrac{\pi}{2}\right]$. Therefore, the inequality (\ref{auxiliaryineq})
is a consequence of the fact that is the decreasing concave curve above the secant line over segment
$\left[0, \nfrac{\pi}{2}\right]$. Based on the inequality (\ref{auxiliaryineq}), we have
$$
\label{startinequality4}
\begin{array}{rcl}
g_{2}(t) > h_{2}(t)
\!\!&\!=\!&\!\!
{\Big (}(\pi^2+\pi-8)\cos^4\!t-\pi{\Big )}\cos t
\left(\nfrac{\pi}{4}
-\nfrac{t}{2}\right)                                                            \\[1.0 ex]
\!\!&\! \!&\!\!
-\,\pi\,\left(\nfrac{\pi}{2}-t\right)^{\!2} + 2\,\pi \cos^2 t,
\end{array}
$$
for $t \in \left( 0, \nfrac{\pi}{2} - 1.1 \right)$.
Thus, we need to prove that
$$
h_{2}(t) > 0,
$$
for $t \in \left( 0, \nfrac{\pi}{2} - 1.1 \right)$.
Let us notice that $h_{2}$ is one mixed trigonometric polynomial
$$
\label{startinequality5}
\begin{array}{rcl}
h_{2}(t)
\!\!&\!=\!&\!\!
\left({\big (}
-\nfrac{\,\pi^2}{2}
-\nfrac{\pi}{2}
+4
{\big )} \, t
+
\nfrac{\,\pi^3}{4}
+
\nfrac{\,\pi^2}{4}
-
2\,\pi
\right) \cos^5 t                                                                \\[1.00 ex]
\!\!&\! \!&\!\!
+
\, 2\,\pi \cos^2 t
+
\left(
\nfrac{\pi}{2} \, t
\!-\!
\nfrac{\,\pi^2}{4}
\right) \cos t
-
\pi \, t^2
+
\pi^2 \, t
-
\nfrac{\,\pi^3}{4}.
\end{array}
$$
For the proof of the inequality $h_{2}(t) > 0$, for $t \in \left( 0, \nfrac{\pi}{2} \!-\! 1.1 \right)$,
we use method from the paper \cite{Malesevic_Makragic_2015}. Using trigonometric multiple angle formulas,
we obtain

\smallskip

\medskip
\noindent
$$
\begin{array}{rcl}
h_{2}(t)
\!&\!=\!&\!
\left(
{\big (}
-\!\nfrac{\pi^2}{32}
-\!\nfrac{\pi}{32}
+\!\nfrac{1}{4}
{\big )} \, t
+\!\nfrac{\pi^3}{64}
+\!\nfrac{\pi^2}{64}
-\!\nfrac{\pi}{8}
\right) \cos 5\,t                                                               \\[1.75 ex]
\!&\! \!&\!
+
\left(
{\big (}
-\!\nfrac{5\,\pi^2}{32}
-\!\nfrac{5\,\pi}{32}
+\!\nfrac{5}{4}
{\big )} \, t
+\!\nfrac{5\,\pi^3}{64}
+\!\nfrac{5\,\pi^2}{64}
-\!\nfrac{5\,\pi}{8}
\right) \cos 3\,t                                                               \\[1.75 ex]
\!&\! \!&\!
+\;
\pi \cos 2\,t
+
\left(
{\big (}
-\!\nfrac{5\,\pi^2}{16}
+\!\nfrac{3\,\pi}{16}
+\!\nfrac{5}{2}
{\big )} \, t
+\!\nfrac{5\,\pi^3}{32}
-\!\nfrac{3\,\pi^2}{32}
-\!\nfrac{5\,\pi}{4}
\right) \cos t                                                                  \\[1.75 ex]
\!&\! \!&\!
-\,\pi t^2
\!+\!\pi^2 t
\!-\!\nfrac{\pi^3}{4}.
\end{array}
$$
For $t \in \left( 0, \nfrac{\pi}{2} - 1.1 \right)$ the following inequalities are true:
$$
\begin{array}{l}
{\big (}
-\!\nfrac{\pi^2}{32}
-\!\nfrac{\pi}{32}
+\!\nfrac{1}{4}
{\big )} \, t
+\!\nfrac{\pi^3}{64}
+\!\nfrac{\pi^2}{64}
-\!\nfrac{\pi}{8}
> 0,                                                                            \\[1.00 ex]
{\big (}
-\!\nfrac{5\,\pi^2}{32}
-\!\nfrac{5\,\pi}{32}
+\!\nfrac{5}{4}
{\big )} \, t
+\!\nfrac{5\,\pi^3}{64}
+\!\nfrac{5\,\pi^2}{64}
-\!\nfrac{5\,\pi}{8}
> 0,                                                                            \\[1.00 ex]
{\big (}
-\!\nfrac{5\,\pi^2}{16}
+\!\nfrac{3\,\pi}{16}
+\!\nfrac{5}{2}
{\big )} \, t
+\!\nfrac{5\,\pi^3}{32}
-\!\nfrac{3\,\pi^2}{32}
-\!\nfrac{5\,\pi}{4}
< 0.
\end{array}
$$
In the purpose of proving that $h_{2}(t) > 0$, for $t \in \left( 0, \nfrac{\pi}{2} \!-\! 1.1 \right)$,
we use the inequalities from \cite{Malesevic_Makragic_2015}:

\bigskip
\noindent
$$
\overline{T}^{\mbox{\footnotesize \rm cos}, 0}_{k}(y)
>
\cos y \,\,{\big (}k \!=\! 0{\big )}
\;\;\;\mbox{and}\;\;\;
\cos y
>
\underline{T}^{\mbox{\footnotesize \rm cos}, 0}_{k}(y)\,\,{\big (}k \!=\! 2{\big )},
$$

\smallskip
\noindent
for $y \in {\big (}\,0,\sqrt{(k+3)(k+4)}\;{\big )}$. Therefore, we get

\break

\noindent
$$
\begin{array}{rcccl}
h_{2}(t)
\!\!&\!>\!&\!\!
P_{3}(t)
\!\!&\!=\!&\!\!
\left(
{\big (}
-\!\nfrac{\pi^2}{32}
-\!\nfrac{\pi}{32}
+\!\nfrac{1}{4}
{\big )} \, t
+\!\nfrac{\pi^3}{64}
+\!\nfrac{\pi^2}{64}
-\!\nfrac{\pi}{8}
\right) \underline{T}^{\mbox{\footnotesize \rm cos}, 0}_{2}(5\,t)               \\[1.2 ex]
\!\!&\! \!&\!\!
\!\!&\! \!&\!\!
+
\left(
{\big (}
-\!\nfrac{5\,\pi^2}{32}
-\!\nfrac{5\,\pi}{32}
+\!\nfrac{5}{4}
{\big )} \, t
+\!\nfrac{5\,\pi^3}{64}
+\!\nfrac{5\,\pi^2}{64}
-\!\nfrac{5\,\pi}{8}
\right) \underline{T}^{\mbox{\footnotesize \rm cos}, 0}_{2}(3\,t)               \\[2.0 ex]
\!\!&\! \!&\!\!
\!\!&\! \!&\!\!
+\;
\pi \, \underline{T}^{\mbox{\footnotesize \rm cos}, 0}_{2}(2\,t)                \\[1.2 ex]
\!\!&\! \!&\!\!
\!\!&\! \!&\!\!
+
\left(
{\big (}
-\!\nfrac{5\,\pi^2}{16}
+\!\nfrac{3\,\pi}{16}
+\!\nfrac{5}{2}
{\big )} \, t
+\!\nfrac{5\,\pi^3}{32}
-\!\nfrac{3\,\pi^2}{32}
-\!\nfrac{5\,\pi}{4}
\right) \overline{T}^{\mbox{\footnotesize \rm cos}, 0}_{0}(t)                   \\[1.2 ex]
\!\!&\! \!&\!\!
\!\!&\! \!&\!\!
-\,\pi t^2
\!+\!\pi^2 t
\!+\!\pi
\!-\!\nfrac{\pi^3}{4},
\end{array}
$$
for $t \in \left( 0, \nfrac{\pi}{2} - 1.1 \right)$.
It is simple to prove that
$$
P_{3}(t)
=
\left(
     \nfrac{35\,\pi^2}{32}
\!+\!\nfrac{35\,\pi}{32}
\!-\!\nfrac{35}{4}
\right) t^3
+
\left(
  -\,\!\nfrac{35\,\pi^3}{64}
\!-\!\nfrac{35\,\pi^2}{64}
\!+\!\nfrac{11\,\pi}{8}
\right) t^2
+
\left(
\nfrac{\pi^2}{2}
+
4
\right) t > 0
$$
for $t \in \left( 0, \nfrac{\pi}{2} - 1.1 \right)$.
Therefore, we conclude that
$$
g_{2}(t) > h_{2}(t) > P_{3}(t) > 0,
$$
for $t \in \left(0, \nfrac{\pi}{2}-1.1\right)$
and consequently that
$$
g(t)>0,
$$
for $t \in \left(1.1,\nfrac{\pi}{2}\right)$, which proves the inequality (\ref{startinequality1}).
The elementary calculus gives
$$
\displaystyle\lim\limits_{x \rightarrow \frac{\pi}{2}_{-}}{
\!\!\nfrac{(\arcsin x / x)^2 + (\arctan x / x ) - 2}{x^3 \, \arctan x}}
\,=\,
\nfrac{\pi^2+\pi-8}{\pi}.
$$
The proof is completed. \stop

\section{Proof of the Conjecture 2}

Let us now observe the inequality (\ref{second inequality}) of the Conjecture \ref{second conjecture} written in the form
\begin{equation}
\label{startinequality2.1}
3+\nfrac{(5\,\pi-12)\,x^3\arctan x}{\pi}-\nfrac{\arctan x}{x}- 2\left(\nfrac{\arctan x}{x}\right)
>
0,
\end{equation}
for $x \in (0,1)$.
Substituting $x=\sin t$ into (\ref{startinequality2.1}), for $t \in \left(0,\nfrac{\pi}{2}\right)$,
we obtain
\begin{equation}
\label{startinequality2.2}
3
+
{\nfrac { \left( ( 5\,\pi -12 ) \sin^{4}  t - \pi \right)
\arctan \left( \sin  t \right) }{\pi \sin t}}
-{\nfrac {2\,t}{\sin t }}
>
0.
\end{equation}
It is enough to prove that
\begin{equation}
\label{startinequality2.3}
\begin{array}{rcl}
g(t)
\!\!&\!=\!&\!\!
3\,\pi \sin t
+((5\pi-12)\sin^{4} t-\pi)\arctan(\sin t)
-2\,\pi \, t
>
0,
\end{array}
\end{equation}
for $t \in \left( 0, \nfrac{\pi}{2} \right)$.
Let us notice that $t = 0$ is zero of the fifth order and $t = \nfrac{\pi}{2}$ is the simple zero of the function $g$.
Furthermore, we differentiate two cases if $t \!\in\! (0, 1.3]$ or $t \!\in\! (1.3, \pi/2)$.

\medskip 
\noindent
{\bf (I)} {\boldmath $t \!\in\! (0, 1.3]$} Based on the inequality (\ref{dobleinequality_4_6}), it may be concluded that
$$
\begin{array}{rcl}
g(t)
\!\!&\!\!>h(t)=\!\!&\!\!
3\,\pi \sin t
+
(5\pi-12) \sin^{4} t \; \underline{T}^{\mbox{\footnotesize \rm arctan}, 0}_{3}(\sin t)    \\[1.0 ex]
\!\!&\! \!&\!\!
-
\pi \; \overline{T}^{\mbox{\footnotesize \rm arctan}, 0}_{5}(\sin t)
-
2\,\pi \, t,
\end{array}
$$
for $t \in \left(0,1.3\right]$. Therefore, we just need to prove

\break

\noindent
$$
\label{h>0}
\begin{array}{rcl}
h(t)
\!\!&\!=\!&\!\!
3\,\pi \sin t
+
(5\pi-12) \sin^{4} \!t
\left(\sin t - \nfrac{1}{3} \sin^3\!t\right )                                   \\[1.5 ex]
\!\!&\! \!&\!\!
-\,
\pi
\left( \sin t - \nfrac{1}{3} \sin^3\!t + \nfrac{1}{5} \sin^5\!t\right )
-
2\,\pi \, t
>
0,
\end{array}
$$
for $t \in \left(0,1.3\right]$. The function $h$ is a mixed trigonometric polynomial function.
For the proof of the inequality $h(t) > 0$, for $t \in \left( 0, 1.3 \right]$,
we use method from the paper \cite{Malesevic_Makragic_2015}.
Using trigonometric multiple angle formulas, we obtain
$$
\begin{array}{rcl}
    h(t)
    \!&\!=\!&\!
       \left( {\nfrac {5\,\pi }{192}} -\nfrac{1}{16} \right) \sin 7\,t
    + \left( {\nfrac {113\,\pi }{960}}-{\nfrac {5}{16}} \right) \sin 5\,t \\[2.0 ex] \!&\! \!&\!
    + \left( -{\nfrac {199\,\pi }{192}}+{\nfrac {39}{16}} \right) \sin 3\,t
    + \left( {\nfrac {833\,\pi }{192}}-{\nfrac {85}{16}} \right) \sin t
     -2\,\pi \,t .
\end{array}
$$
We also need inequalities from the paper \cite{Malesevic_Makragic_2015}:
$$
\overline{T}^{\mbox{\footnotesize \rm sin}, 0}_{k}(y)
>
\sin y\,\,{\big (}k\!=\!9{\big )}
\;\;\;\mbox{and}\;\;\;
\sin y
>
\underline{T}^{\mbox{\footnotesize \rm sin}, 0}_{k}(y)\,\,{\big (}k\!=\!7, 15, 19{\big )},
$$
for $y \in {\big (}\,0,\sqrt{(k+3)(k+4)}\;{\big )}$.
Putting things together, we get
$$
\begin{array}{rcccl}
h(t)
\!&\!>\!&\!
P_{19}(t)
\!&\!=\!&\!
       \left( {\nfrac {5\,\pi }{192}} -\nfrac{1}{16} \right) \; \underline{T}^{\mbox{\footnotesize \rm sin}, 0}_{19}(7 \; t)   \\[1.5 ex]
&&&& + \left( {\nfrac {113\,\pi }{960}}-{\nfrac {5}{16}} \right) \underline{T}^{\mbox{\footnotesize \rm sin}, 0}_{15}(5 \; t)  \\[1.5 ex]
&&&& + \left( -{\nfrac {199\,\pi }{192}}+{\nfrac {39}{16}} \right) \overline{T}^{\mbox{\footnotesize \rm sin}, 0}_{9}(3 \; t)  \\[1.5 ex]
&&&& + \left( {\nfrac {833\,\pi }{192}}-{\nfrac {85}{16}} \right) \underline{T}^{\mbox{\footnotesize \rm sin}, 0}_{7}(t)
       -2\,\pi \,t,
\end{array}
$$

\smallskip
\noindent
for $t \in \left(0,1.3\right]$.
Hence, we only have to prove that
$$
\begin{array}{rcl}
P_{19}(t)
\!&\!=\!&\!
  \left( {-{\nfrac {232630513987207\,\pi }{95330037871411200}} + \nfrac {232630513987207}{39720849113088000}} \right) {t}^{19} \\[1.75 ex] \!&\!\!&\!
+ \left( {\nfrac {4747561509943\,\pi }{278742800793600}} - {\nfrac {4747561509943}{116142833664000}} \right) {t}^{17} \\[1.75 ex] \!&\!\!&\!
+ \left( {-\nfrac {111034112797\,\pi }{1141243084800}} + {\nfrac {612518675071}{2615348736000}} \right) {t}^{15} \\[1.75 ex] \!&\!\!&\!
+ \left( {\nfrac {25601647133\,\pi }{59779399680}} - {\nfrac {585184807}{566092800}} \right) {t}^{13} \\[1.75 ex] \!&\!\!&\!
+ \left( {-\nfrac {549507467\,\pi }{383201280}} + {\nfrac {277683421}{79833600}} \right) {t}^{11} \\[1.75 ex] \!&\!\!&\!
+ \left( {\nfrac {34570249\,\pi }{9953280}} - {\nfrac {9870319}{1161216}} \right) {t}^{9} \\[1.75 ex] \!&\!\!&\!
+ \left( {-\nfrac {473\,\pi }{84}} + \mbox{\small 14} \right) {t}^{7} \\[1.75 ex] \!&\!\!&\!
+ \left( {\nfrac {93\,\pi }{20}} - \mbox{\small 12} \right) {t}^{5}
>
0,
\end{array}
$$
for $t \in \left(0,1.3\right]$. Let us notice that $P_{19}(t) = \nfrac{t^5}{23355859278495744000} P_{14}(t)$,
where
$$
\begin{array}{rcl}
P_{14}(t)
\!\!&\!\!=\!\!&\!\!
\big(\! -\!56994475926865715\,\pi +136786742224477716 \,\!\big) {t}^{14} \\[1.0 ex] \!\!&\!\! \!\!&\!\!
+ \, \big(\!\, 397798178918123970\,\pi -954715629403497528 \,\!\big) {t}^{12} \\[1.0 ex] \!\!&\!\! \!\!&\!\!
+ \, \big(\! -\!2272344207942188160\,\pi +5469977974061251584 \,\!\big) {t}^{10} \\[1.0 ex] \!\!&\!\! \!\!&\!\!
+ \, \big(\!\, 10002584016180806400\,\pi -24143557388833935360 \,\!\big) {t}^{8} \\[1.0 ex] \!\!&\!\! \!\!&\!\!
+ \, \big(\! -\!33492109086208281600\,\pi +81238161686899875840 \,\!\big) {t}^{6} \\[1.0 ex] \!\!&\!\! \!\!&\!\!
+ \, \big(\!\, 81120783386638195200\,\pi -198524461941501696000 \,\!\big) {t}^{4} \\[1.0 ex] \!\!&\!\! \!\!&\!\!
+ \, \big(\! -\!131515731413434368000\,\pi +326982029898940416000 \,\!\big) {t}^{2} \\[1.0 ex] \!\!&\!\! \!\!&\!\!
+ \, 108604745645005209600\,\pi -280270311341948928000.
\end{array}
$$
Let us introduce substitution $z\!=\!t^2$, for $z\in (0, 1.69]$ and prove $P_{7}(z) = P_{14}(\sqrt{z}\,) > 0$.
It is enough to observe that the third\,-\,order derivative polynomial
$P_{7}'''$ does not have real roots according to the Ferrari's formula and $P_{7}'''(0) < 0$
yields $P_{7}'''(z) < 0$ for every $z \in (0, 1.69]$. Thus
the second\,-\,order derivative polynomial $P_{7}''$ is strictly decreasing function with unique
real root  $z_1 = 1.834 \ldots > 1.69$. From $P_{7}''(1.69) > 0$ follows that $P_{7}''(z) > 0$ for
$z \in (0, 1.69]$, thus derivative polynomial $P_{7}'$ is strictly increasing function for
$z \in (0, 1.69]$. As $P_{7}'(1.69) < 0$ then $P_{7}'(z) < 0$ for every $z \in (0, 1.69]$.
From $P_{7}'(z) < 0$ follows that $P_{7}$ is strictly decreasing function (with real
root $z_2 = 1.870 \ldots > 1.69$). As $P_{7}(1.69) > 0$, we conclude $P_{7}(z) > 0$
for $z \in (0, 1.69]$, i.e. $P_{14}(t) > 0$ for $t \in (0,1.3]$. Finally, from
$$
g(t) > h(t) > P_{19}(t) = \nfrac{\,t^5}{23355859278495744000} P_{14}(t)
$$
follows that $g(t) > 0$ for $t \in (0,1.3]$.

\bigskip  
\noindent
{\bf (II)} {\boldmath $t \!\in\! \left(1.3, \pi\!\,/\,\!2\right)$}
We transform the inequality (\ref{startinequality2.3}), for $t \in \left(1.3,\nfrac{\pi}{2}\right)$,
to the~\mbox{inequality}
$$
\label{startinequality3}
\begin{array}{rcl}
g_{2}(t)
=
g(\nfrac{\pi}{2} - t)
\!\!&\!=\!&\!\!
\left( ( 5\,\pi -12 ) \cos^{4} t -\pi  \right) \arctan \left( \cos t \right)    \\[1.5 ex]
\!\!&\! \!&\!\!
+\,3\,\pi \cos t \, - 2 \, \pi \, \left( \nfrac{\pi}{2}-t \right)
>
0,
\end{array}
$$
for $t \in \left( 0, \nfrac{\pi}{2} - 1.3 \right) = {\big (} 0, 0.270 \ldots {\big )}$. Let us notice that
$\left( 5\,\pi -12 \right) \cos^{4} t -\pi > 0$ is true for $t \in \left( 0, \nfrac{\pi}{2} - 1.3 \right)$.
Based on the inequality (\ref{auxiliaryineq}) we have
$$
\label{startinequality4}
\begin{array}{rcl}
g_{2}(t) > h_{2}(t)
\!\!&\!\!=\!\!&\!\!
\left( ( 5\,\pi - 12 ) \cos^{4} \; t -\pi  \right)  \left( \nfrac{\pi} {4}-\nfrac{t}{2} \right)   \\[1.5 ex]
\!\!&\! \!&\!\!
+\,3\,\pi \cos t \, - 2 \, \pi \, \left( \nfrac{\pi}{2}-t \right),
\end{array}
$$
for $t \in \left( 0, \nfrac{\pi}{2} - 1.3 \right)$; so it should be proved
$$
\label{h2>0}
h_{2}(t) > 0,
$$
for $t \in \left( 0, \nfrac{\pi}{2} - 1.3 \right)$.
Notice that $h_{2}$ is one mixed trigonometric polynomial
$$
\label{startinequality5}
\begin{array}{rcl}
h_{2}(t)
\!\!&\!=\!&\!\!
\left( {\big (}6-\nfrac{5\,\pi}{2}{\big )}\,t + \nfrac{5\,{\pi }^{2}}{4}-3\,\pi \right)\cos^{4}t  \\[1.5 ex]
\!\!&\! \!&\!\!
+\,3\,\pi \cos t + \nfrac{5\,\pi}{2} \,t - \nfrac{5\,{\pi }^{2}}{4}.
\end{array}
$$
For the proof of the inequality $h_{2}(t) > 0$, for $t \in \left( 1.3, \pi/2 \right)$, we use method from
the paper \cite{Malesevic_Makragic_2015}. Using trigonometric multiple angle formulas, we obtain
$$
\begin{array}{rcl}
h_{2}(t)
\!&\!=\!&\!
  \left( {\big (}-\nfrac {5\,\pi }{16} + \nfrac{3}{4}{\big )} \, t + {\nfrac {5\,{\pi }^{2}}{32}} - \nfrac{3}{8}\,\pi  \right) \cos 4\,t \\[1.2 ex] \!&\! \!&\!
+ \left( {\big (} -\nfrac{5\,\pi}{4} + 3{\big )} \, t + \nfrac{5\,{\pi }^{2}}{8} - \nfrac{3\,\pi}{2} \right) \cos 2\,t                  \\[2.0 ex] \!&\! \!&\!
+ \,3\,\pi \cos t                                                                                                          \\[1.2 ex] \!&\! \!&\!
+ {\big (} \nfrac{25\,\pi }{16} + \nfrac{9}{4} {\big )} \, t
- \nfrac {25\,{\pi }^{2}}{32} - {\nfrac {9\,\pi }{8}}.
\end{array}
$$
For $t \in \left( 0, \nfrac{\pi}{2} - 1.3 \right)$ we have
$$
\begin{array}{l}
{\big (} -\nfrac {5\,\pi}{16} + \nfrac{3}{4} {\big )} \, t + \nfrac {5\,{\pi }^{2}}{32} - \nfrac{3\,\pi}{8} > 0,                         \\[1.5 ex]
{\big (} -\nfrac{5\,\pi}{4} +3 {\big )} \, t + \nfrac{5\,{\pi }^{2}}{8} - \nfrac{3\,\pi}{2} > 0,                                         \\[1.5 ex]
\end{array}
$$
Having in mind inequality from \cite{Malesevic_Makragic_2015}:

\smallskip
$$
\cos y
>
\underline{T}^{\mbox{\footnotesize \rm cos}, 0}_{2}(y)\,,
$$

\smallskip
\noindent
for $y \in {\big (}\,0,\sqrt{30}\;{\big )}$,
we conclude that
$$
\begin{array}{rcccl}
h_{2}(t)
\!\!&\!>\!&\!\!
P_{3}(t)
\!\!&\!=\!&\!\!
  \left( {\big (} -\nfrac {5\,\pi }{16} + \nfrac{3}{4} {\big )} \, t + \nfrac {5\,{\pi }^{2}}{32} - \nfrac{3\,\pi}{8} \right) \underline{T}^{\mbox{\footnotesize \rm cos}, 0}_{2}(4\;t)  \\[1.2 ex] \!\!&\! \!&\!\! \!\!&\! \!&\!\!
+ \left( {\big (} -\nfrac{5\,\pi}{4} +3 {\big )} \, t + \nfrac{5\,{\pi }^{2}}{8} - \nfrac{3\,\pi}{2} \right) \underline{T}^{\mbox{\footnotesize \rm cos}, 0}_{2}(2\;t)                   \\[2.0 ex] \!\!&\! \!&\!\! \!\!&\! \!&\!\!
+ \,3\,\pi\;\underline{T}^{\mbox{\footnotesize \rm cos}, 0}_{2}(t)                                                                                                                       \\[1.2 ex] \!\!&\! \!&\!\! \!\!&\! \!&\!\!
+ {\big (} \nfrac{25\,\pi}{16} + \nfrac{9}{4} {\big )} \, t - {\nfrac {25\,{\pi }^{2}}{32}}-{\nfrac {9\,\pi }{8}},
\end{array}
$$
for $t \in \left( 0, \nfrac{\pi}{2} - 1.3 \right)$. It is simple to prove that
$$
P_{3}(t)
=
\left( 5\,\pi -12 \right) {t}^{3} + {\big (} -\nfrac{5\,{\pi }^{2}}{2}+\nfrac{9\,\pi}{2} {\big )}\,t^{2} + 6\,t > 0,
$$
for $t \in \left( 0, \nfrac{\pi}{2} - 1.3 \right)$.
Therefore, we conclude that
$$
g_{2}(t) > h_{2}(t) > P_{3}(t) > 0,
$$
for $t \in \left(0, \nfrac{\pi}{2}-1.3\right)$
and consequently that
$$
g(t)>0,
$$
for $t \in \left(1.3,\nfrac{\pi}{2}\right)$, which proves the inequality (\ref{startinequality2.2}).
The elementary calculus proposes
$$
\displaystyle\lim\limits_{x \rightarrow \frac{\pi}{2}_{-}}{
\!\!\nfrac{2(\arcsin x / x) + (\arctan x / x ) - 3}{x^3 \arctan x}}
\,=\,
\nfrac{5\pi-12}{\pi}.
$$

\smallskip
\noindent
Therefore, the proof of the second conjecture is also completed. \stop

\section{Results, Discussion and Conclusions}

Let us emphasize that the method from \cite{Malesevic_Makragic_2015} was here applied for proving two conjectures stated by C.-P. Chen \cite{Chao-Ping Chen}.
In the paper \cite{Malesevic_Makragic_2015} the open problem stated by Z.-J. Sun and L. Zhu \cite{Sun_Zhu_2011} was also proved in the same manner.
We expect that the method will be useful in solving some others problem concerning inequalities
which can be reduced to some mixed trigonometric inequalities.



\medskip

\medskip
\noindent
\textbf{Acknowledgements.}
The research is partially supported by the Ministry of Education and Science, Serbia, Grants No.
174032 and 44006.

\break

\break

\end{document}